\documentclass[a4paper]{article}
\usepackage{amsfonts}
\usepackage{amssymb}
\usepackage{amsmath}
\usepackage{latexsym}
\usepackage{calrsfs}
\newtheorem{theorem}{Theorem}[section]
\newtheorem{proposition}[theorem]{Proposition}
\newtheorem{corollary}[theorem]{Corollary}

\newtheorem{example}[theorem]{Example}
\newtheorem{definition}[theorem]{Definition}
\newtheorem{remark}[theorem]{Remark}
\newtheorem{problem}{Problem}

\newcommand{\proof}{\medskip \noindent {\bf Proof. \ \ }}
\newcommand{\qed}{\null\hfill $\Box\;\;$ \medskip}

\begin{document}

\parbox{1mm}

\begin{center}
{\bf {\sc \Large On distance distribution functions-valued
submeasures related to aggregation functions}}
\end{center}

\vskip 12pt

\begin{center}
{\bf Lenka HAL\v{C}INOV\'A, Ondrej HUTN\'IK and Radko
MESIAR}
\end{center}


\hspace{5mm}\parbox[t]{10cm}{\fontsize{9pt}{0.1in}\selectfont\noindent{\bf
Abstract.} Probabilistic submeasures generalizing the classical
(numerical) submeasures are introduced and discussed in connection
with some classes of aggregation functions. A special attention is
paid to triangular norm-based probabilistic submeasures and
semi-copula-based probabilistic submeasures. Some algebraic
properties of classes of such submeasures are also studied. }
\vskip 24pt

\section{Introduction and motivations}

In recent years non-additive set functions have attracted much
attention in pure mathematics as well as in various applications.
This class includes well-known set functions such as submeasures,
Dobrakov submeasures and semi-measures, fuzzy measures, null
additive set functions, etc. As a larger overview of non-additive
set functions we recommend monographs~\cite{Pap} and~\cite{WK}, or
several chapters in the handbook~\cite{PapH}.

The study of submeasures was initiated in the second half of the
last century by Orlicz and developed by
Drewnowski~\cite{Drewnowski} from the topological point of view.
In fact, many classical objects of measure theory, as e.g.
variations and semi-variations of vector measures, are
submeasures. This classical object of measure theory is defined as
follows. Let $\Sigma$ be a ring of subsets of a fixed (non-empty)
set $\Omega$ and $\overline{\mathbb{R}}_+=[0,+\infty]$ be the
extended non-negative real half-line. A mapping $\eta: \Sigma \to
\overline{\mathbb{R}}_+$ satisfying the conditions
\begin{itemize}
\item[(i)] $\eta(\emptyset) = 0$; \item[(ii)] $\eta(E)\leq
\eta(F)$ for $E,F \in \Sigma$ such that $E \subset F$;
\item[(iii)] $\eta(E \cup F) \leq \eta(E) + \eta(F)$ whenever $E,F
\in \Sigma$.
\end{itemize} \noindent is said to be a \textit{numerical submeasure}
on $\Sigma$.

In our previous papers~\cite{HutMes} and~\cite{HalHutMes} we have
introduced and investigated submeasure notions related to
probabilistic metric spaces (PM-spaces, for short), see~\cite{SS}.
Our considerations of a submeasure notion in paper~\cite{HutMes}
were closely related to the Menger PM-space $(\Omega, \mathcal{F},
\tau_T)$ where $\tau_T$ is the triangle function in the form
\begin{equation}\label{tauT}
\tau_T(G,H)(x) = \sup_{u+v=x}T(G(u),H(v))
\end{equation}\noindent with $T$ being a left-continuous t-norm. The
associated submeasure notion was defined as follows,
see~\cite[Definition~3]{HutMes} and Section~\ref{sectionprem} for
necessary notations.

\begin{definition}\label{defT-sub}\rm
Let $T: [0,1]^2 \to [0,1]$ be a t-norm, and $\Sigma$ a ring of
subsets of $\Omega\neq\emptyset$. A mapping $\gamma: \Sigma \to
\Delta^+$ (where $\gamma(E)$ is denoted by $\gamma_E$) such that
\begin{itemize}
\item[(a)] if $E=\emptyset$, then $\gamma_E(x) =
\varepsilon_0(x)$, $x>0$; \item[(b)] if $E \subset F$, then
$\gamma_E(x) \geq \gamma_F(x)$, $x>0$; \item[(c)] $\gamma_{E\cup
F}(x+y) \geq T(\gamma_E(x), \gamma_F(y))$, \,\,$x,y>0$, $E,F \in
\Sigma$,
\end{itemize}\noindent is said to be a \textit{$\tau_T$\,-submeasure}.
\end{definition}

From this definition is obvious that the probabilistic submeasure
is a certain (non-additive) set function taking values in the set
of distribution functions of non-negative random variables. The
attribute "submeasure" reflects the property~(c) which is a
"probabilistic" version of the classical subadditivity. The origin
of this notion comes from the fact that it works in such
situations in which we have only a \textit{probabilistic
information} about measure of a set (recall a similar situation in
the framework of information measures). 
For example, if rounding of reals is considered, then the uniform
distributions over intervals describe our information about the
measure of a set. 

On the other hand, $\tau_T$-submeasures can be seen as fuzzy
number-valued submeasures. In this case the value $\gamma_E$ can
be seen as a non-negative $LT$-fuzzy number, see~\cite{DKMP},
where $\tau_T(\gamma_E,\gamma_F)$ corresponds to the $T$-sum of
fuzzy numbers $\gamma_E$ and $\gamma_F$. Moreover, each
$\tau_T$-submeasure $\gamma$ with the minimum t-norm $T=M$
(in~\cite{HutMes} we call it \textit{universal
$\tau_T$-submeasure}) can be represented by means of a
non-decreasing system $(\eta_\alpha)_{\alpha \in [0,1]}$ of
numerical submeasures (compare the horizontal representation
$(S_\alpha)_{\alpha \in [0,1]}$ of a fuzzy subset $S$), where
$$\gamma_E(x) = \sup\{\alpha \in [0,1];\, \eta_\alpha(E)\leq x\}, \quad E\in\Sigma.$$

\begin{example}\rm
Let $\eta$ be a numerical submeasure on $\Sigma$. Then for each
$E\in\Sigma$ the mapping
$$\gamma_{E}(x) =
1-\exp\left(-\left[\frac{x}{\lambda\,\eta(E)}\right]^{k}\right),
\quad x>0, \lambda>0, k>0,$$ corresponds to a cumulative
distribution function of the Weibull distribution $W(\lambda,k)$
with parameters $\lambda, k$. Especially, for $k=1$ we get the
(universal) $\tau_{T}$-submeasure corresponding to a distribution
function of exponential distribution $E(\lambda)$ with parameter
$\lambda$. Note that the standard conventions for the arithmetic
operations on $\overline{\mathbb{R}}_+$ are considered, such as
$0\cdot (+\infty) = 0/0 = 0$.
\end{example}

Naturally, we may ask about possibility to extend our
considerations from Menger PM-spaces to wider spaces with
different triangular functions instead of~(\ref{tauT}). Note that
similar considerations were introduced and discussed in the
framework of probabilistic metric spaces, see for example the
monograph~\cite{GChR}. For such reasons in paper~\cite{HalHutMes}
we have provided a generalization of $\tau_{T}$-sub\-measu\-res
which involves suitable operations $L$ replacing the standard
addition $+$ on $\overline{\mathbb{R}}_+$ such that the underlying
function~(\ref{tauT}) is a triangle function and thus the
underlying space is the so-called $L$-Menger PM-space. Since
t-norms are rather special operations on the unit interval $[0,
1]$, we have also mentioned few possible generalizations of a
submeasure notion based on aggregation operators and convolution
of distance distribution functions, i.e., such submeasures which
can be used in non-Menger PM-spaces (e.g., in the Wald spaces),
but also in wider class of PM-spaces.

The aim of this paper is a further generalization of the concept
of probabilistic submeasures. In particular, triangular norms
applied in (c) of Definition~\ref{defT-sub} are used as binary
functions only, and thus their associativity is a superfluous
constraint. Therefore, a more general aggregation function can be
used here (compare, e.g., the case of fuzzy logics, where the
triangular norms can be replaced by (quasi-)copulas as discussed
in~\cite{HM}).

The paper is organized as follows: in Section~\ref{sectionprem} we
recall some basic and necessary notions which will be used in this
paper. Then in Section~\ref{sectiont-norms} we investigate further
properties of triangular norm-based probabilistic submeasures
which generalize some results obtained in our previous papers.
Passing from triangular norms to their natural
extension/modification in the form of copulas, quasi-copulas and
semi-copulas we study in Section~\ref{sectioncopulas} notion of
submeasures related to these aggregation functions. In the whole
paper a number of examples is presented. The lattice structure of
spaces of semi-copula and quasi-copula-based submeasures is also
discussed.

\section{Basic notions and definitions}\label{sectionprem}

The class of all distance distribution functions (distribution
functions of non-negative random variables) will be denoted by
$\Delta^+$. A \textit{triangle function} is a function $\tau:
\Delta^+ \times \Delta^+ \to \Delta^+$ which is symmetric,
associative, non-decreasing in each variable and has
$\varepsilon_0$ as the identity, where $\varepsilon_0$ is the
distribution function of Dirac random variable concentrated in
point 0. More precisely, for $a\in[0,+\infty[$ we put
$$\varepsilon_a(x) =
\begin{cases}
1& \textrm{for}\,\, x>a, \\
0 & \textrm{otherwise}.
\end{cases}$$ Clearly, $(\Delta^+, \tau)$ is an
Abelian semigroup with the identity $\varepsilon_0$. A
\textit{triangular norm}, shortly a t-norm, is a commutative
lattice ordered semi-group on $[0,1]$ with identity 1. The most
important t-norms are the minimum $M(x,y) = \min\{x,y\}$, the
product $\Pi(x,y) = xy$, the {\L}ukasiewicz $W(x,y) =
\max\{x+y-1,0\}$ and the drastic product $$D(x,y) =
\begin{cases}
\min\{x,y\} & \textrm{for}\,\, \max\{x,y\}=1 \\
0 & \textrm{otherwise}.
\end{cases}$$ For more information about t-norms we refer the book~\cite{KMP}.
We denote by $\mathcal{T}$ the class of all t-norms.

Triangular norms are a rather special case of aggregation
functions on $[0,1]$. A binary \textit{aggregation function} $A:
[0,1]^2\to [0,1]$ is a non-decreasing function in both components
with the boundary conditions $A(0,0)=0$ and $A(1,1)=1$. The class
of all binary aggregation functions will be denoted by
$\mathcal{A}$. For more details on aggregation functions we
recommend a recent monograph~\cite{Grabish}.

Let us denote by $\mathcal{L}$ the set of binary operations on
$\overline{\mathbb{R}}_+$ such that

\begin{itemize} \item[(i)] $L$ is commutative and associative; \item[(ii)] $L$ is jointly strictly
increasing, i.e., for all $u_1, u_2, v_1,
v_2\in\overline{\mathbb{R}}_+$ with $u_1<u_2$, $v_1<v_2$ holds
$L(u_1,v_1)<L(u_2,v_2)$; \item[(iii)] $L$ is continuous on
$\overline{\mathbb{R}}_+\times\overline{\mathbb{R}}_+$;
\item[(iv)] $L$ has $0$ as its neutral element.
\end{itemize}\noindent Observe that $L\in\mathcal{L}$
is a jointly increasing pseudo-addition on
$\overline{\mathbb{R}}_+$ in the sense of~\cite{SM}. The usual
examples of operations in $\mathcal{L}$ are
\begin{align*}
K_\alpha(x,y) & = (x^\alpha+y^\alpha)^{\frac{1}{\alpha}}, \quad
\alpha>0, \\ K_\infty(x,y) & = \max\{x,y\}.
\end{align*}Note that although $\max\{x,y\}\in\mathcal{L}$, its "counterpart"
$\min\{x,y\}$ is not a member of $\mathcal{L}$, because
$\min\{x,y\}$ does not have $0$ as its neutral element. With
$(L,A)\in\mathcal{L}\times\mathcal{A}$ the general form
of~(\ref{tauT}) is as follows $$\tau_{L,A}(G,H)(x) =
\sup_{L(u,v)=x} A(G(u), H(v)).$$ Note that the left-continuity of
$A$ ensures that $\tau_{L,A}$ is a binary operation on $\Delta^+$.
However, $\tau_{L,A}$ need not be associative in general, but it
has good properties on $\Delta^+$.

Now we introduce the following probabilistic submeasure notion in
its general form (note that neither the left-continuity of a
t-norm $T$ nor of an aggregation function $A$ is required in what
follows).

\begin{definition}\rm\label{defLT-sub}
Let $(L,A)\in\mathcal{L}\times\mathcal{A}$ and $\Sigma$ be a ring
of subsets of $\Omega\neq \emptyset$. A mapping $\gamma: \Sigma
\to \Delta^+$ such that
\begin{itemize}
\item[(a')] $\gamma_E(x) = \varepsilon_0(x)$, $x>0$; \item[(b')]
$\gamma_E(x) \geq \gamma_F(x)$, $x>0$ whenever $E\subset F$;
\item[(c')] $\gamma_{E\cup F}(L(x,y)) \geq A(\gamma_E(x),
\gamma_F(y))$, \,\,$x,y>0$, $E,F \in \Sigma$,
\end{itemize}\noindent is said to be a \textit{$\tau_{L,A}$\,-submeasure}.
\end{definition}

\noindent If $L=K_1$, then its index is usually omitted, and we
simply speak about $\tau_A$-submeasure. Clearly, for $A=T$ (a
left-continuous t-norm), and $L=K_1$ the $\tau_{L,A}$-submeasure
reduces to $\tau_T$-submeasure from~\cite{HutMes}. For instance,
for $L=K_\infty$ we get a $\tau_{\max, T}$-submeasure related to a
non-Archimedean Menger PM-space $(\Omega, \mathcal{F}, \tau_{\max,
T})$. It is worth to note that in this case condition~(c') reads
as follows
$$\gamma_{E\cup F}(s) \geq T(\gamma_E(s), \gamma_F(s)), \quad s >0,
E,F \in \Sigma.$$

\begin{remark}\rm
In general, $L\in\mathcal{L}$ if and only if there is a (possibly
empty) system $(]a_k,b_k[)_{k \in K}$ of pairwise disjoint open
subintervals of $]0,+\infty[$, and a system $(\ell_k)_{k \in K}$
of increasing bijections $\ell_k: [a_k,b_k]\to
\overline{\mathbb{R}}_+$ so that
$$L(x,y) =
\begin{cases}
\ell_{k}^{-1}(\ell_k(x)+\ell_k(y)) & \mbox{if } (x,y) \in ]a_k,b_k[^2,\\
\max\{x,y\} & \mbox{otherwise}.
\end{cases}$$ For more details see~\cite{KMP}. For $L=K_\alpha\in\mathcal{L}$ and $A=T\in\mathcal{T}$ we have
$$\tau_{K_\alpha,T} (G,H)(x) = \tau_T(G,H)(x^\alpha),$$ which motivates us to say that for
$L\in\mathcal{L}$ generated by a strictly increasing bijection
$\ell: \overline{\mathbb{R}}_+\to \overline{\mathbb{R}}_+$, we
denote $L = K_\ell$, we have
$$\tau_{K_\ell,T} (G,H)(x) = \tau_T(G,H)(\ell(x)).$$ In this light we
have
\begin{displaymath}
\gamma_{E\cup F}(L(x,y))=\gamma_{F\cup E}(L(y,x)) \geq
\max\{A(\gamma_E(x), \gamma_F(y)), A(\gamma_F(y), \gamma_E(x))\},
\end{displaymath} for $x,y>0$, $E,F\in\Sigma$. So, we may (equivalently) take the
symmetrization $$A_{\rm{sym}}(u,v)=\max\{A(u,v), A(v,u)\}$$
instead of $A\in\mathcal{A}$.
\end{remark}

Easily, by standard methods of measure theory it is possible to
extend a $\tau_{L,A}$-submeasure $\gamma$ from a ring
$\Sigma\subset \mathfrak{P}(\Omega)$ of subsets of $\Omega\neq
\emptyset$ to a set function $\gamma^*:
\mathfrak{P}(\Omega)\to\Delta^+$ as follows
$$\gamma_E^*(x) = \sup\{\gamma_F(x); \,\,E\subseteq F\in\Sigma\}, \quad x>0, E\in\Omega.$$

\begin{problem}
Fix $(L,A)\in\mathcal{L}\times\mathcal{A}$ and let $\gamma$ be a
$\tau_{L,A}$-submeasure on a ring $\Sigma$ of subsets of
$\Omega\neq \emptyset$. Is the (Jordan) extension $\gamma^*$ of
$\gamma$ also a $\tau_{L,A}$-submeasure on $\mathfrak{P}(\Omega)$?
\end{problem}

For better readability in what follows we use the following
convention: since $\Delta^+$ is the set of all distribution
functions with support $\overline{\mathbb{R}}_+$, we state the
expression for a $\tau_{L,A}$\,-submeasure $\gamma: \Sigma \to
\Delta^+$ with $(L,A)\in\mathcal{L}\times\mathcal{A}$ just for
positive values of $x$. Usually we also omit the information "for
$x>0$" if it is not necessary and in accordance with our
convention. In case $x\leq 0$ we always suppose
$\gamma_{\cdot}(x)=0$.

Furthermore, in the whole paper $\Theta_{L,A}$ denotes the set of
all $\tau_{L,A}$-submeasures on $\Sigma$ for a fixed
$(L,A)\in\mathcal{L}\times \mathcal{A}$ and
$$\Theta_{\mathcal{L},\mathcal{A}} = \{\Theta_{L,A};\,
(L,A)\in\mathcal{L}\times\mathcal{A}\}$$ the set of all
$\tau_{L,A}$-submeasures on $\Sigma$ for all possible pairs
$(L,A)\in\mathcal{L}\times \mathcal{A}$ (or, the "superset" of all
sets of $\tau_{L,A}$-submeasures on $\Sigma$). Here also, as a
convention, we omit the index $L=K_1$ and write $\Theta_A$ instead
of $\Theta_{K_1,A}$.

\begin{example}\rm
For the set $\Omega = \{\omega_1,\omega_2\}$ and positive
constants $a,b,c$ such that $c \leq \min\{a,b\}$, put
\begin{align*}
\gamma_{\omega_1}(x) & = \max\{0,1 - e^{-ax}\}, \\
\gamma_{\omega_2}(x) & = \max\{0,1 - e^{-bx}\}, \\
\gamma_\Omega(x) & = \max\{0,1 - e^{-cx}\}.
\end{align*} Then $\gamma \in \Theta_A$, $A\in\mathcal{A}$, if and only if $$A(u,v) \leq
(1-u)^\frac{c}{a} \cdot (1-v)^\frac{c}{b}.$$ Hence, for $a = b =
\frac{3}{2}c$, we have $\gamma \in \Theta_\Pi$, but $\gamma \notin
\Theta_M$. Observe also that $\gamma$ is not related to any
numerical submeasure, see~\cite{HutMes}.
\end{example}

\begin{example}\rm
For a positive real number $p$ consider the class
$\mathbf{M}_p\subset\mathcal{A}$ which is usually called the
\textit{$p$-mean} (or, the H\"{o}lder mean) and is defined as
$$\mathbf{M}_p(x,y)=\left(\frac{x^p+y^p}{2}\right)^{1/p}, \quad x,y\geq 0.$$ If $\eta$ is a numerical submeasure
on $\Sigma$, then $\gamma\in\Theta_{\mathbf{M}_p}$, where
$$\gamma_{E}(x)=2^{-1/p}\Biggl(1+\left(\max\left\{\min\left\{\sqrt[p]{\max\left\{1+p(x-\eta(E)),0\right\}},\, 1\right\}, 0\right\}\right)^p\Biggr)^{1/p}
$$ and $E\in\Sigma$. Since $\lim\limits_{p\to 0}\mathbf{M}_p = \mathbf{G}$, the
\textit{geometric mean}, then $\gamma\in\Theta_{\mathbf{G}}$ has
the form
$$\gamma_{E}(x) = \sqrt{\min\left\{e^{x-\eta(E)}, 1\right\}}, \quad
E\in\Sigma.$$ Also, for $p=1$, resp. $p=2$, the $p$-mean is
nothing but the \textit{arithmetic mean} $\mathbf{A}$, resp. the
\textit{quadratic mean} $\mathbf{Q}$, and therefore we easily get
the corresponding $\tau_{\mathbf{A}}$-, resp.
$\tau_{\mathbf{Q}}$-submeasure.
\end{example}

\section{Triangular norm-based submeasures}\label{sectiont-norms}

In what follows we use the usual point-wise order $\leq$ between
real-valued functions. Since $\gamma\in \Delta^+$ is
non-decreasing, then for a fixed $T\in\mathcal{T}$ each
$\tau_{L_1,T}$-submeasure is a $\tau_{L_2, T}$-submeasure whenever
$L_1\leq L_2$. Moreover, if $T_2\leq T_1$ (it is usually said that
$T_2$ is a \textit{weaker} t-norm than $T_1$, or $T_1$ is
\textit{stronger} than $T_2$, see~\cite{KMP}), then each
$\tau_{L_1,T_1}$-submeasure is a $\tau_{L_2,T_2}$-submeasure.
In accordance with this motivation introduce the order $\ll$ on
$\Theta_{\mathcal{L},\mathcal{T}}$ as follows
$$\Theta_{L_1, T_1} \ll \Theta_{L_2, T_2} \quad \textrm{if and
only if} \quad L_1 \leq L_2\,\,\textrm{and}\,\, T_2\leq T_1.$$
Then $(\Theta_{\mathcal{L},\mathcal{T}}, \ll)$ is a partially
ordered set and for each $(L,T)\in\mathcal{L}\times\mathcal{T}$ we
have
$$\Theta_{L,M} \ll \Theta_{L,T} \ll \Theta_{L,D}.$$ Note that
for $L=K_1$ the order $\ll$ on $\Theta_{\mathcal{T}}$ is nothing
but order-inverted image of the point-wise order $\leq$ of
t-norms.

\begin{table}
\begin{center}{\small
\begin{tabular}{|l||l|}
\hline & \\ \textbf{Family of t-norms} & \textbf{Corresponding
family of $\tau_{T}$-submeasures} \\ & \\
\hline

\hline & \\ \textit{Acz\'el-Alsina t-norms} & $\gamma_E^{AA,0}(x)
= \varepsilon_{\eta(E)}(x)$ \\ $T_{\lambda}^{AA}$, $\lambda \in
[0,+\infty[$ & $\gamma_E^{AA,\lambda}(x) =
\exp\left(-\Bigl[\max\{\eta(E)-x,0\}\Bigr]^{1/\lambda}\right)$ \\
& \\

\hline & \\ \textit{Dombi t-norms} & $\gamma_E^{D,0}(x) =
\gamma_E^{AA,0}(x)$ \\
$T_{\lambda}^{D}$, $\lambda \in [0,+\infty[$ &
$\gamma_E^{D,\lambda}(x)
=\left(1+\Bigl[\max\{\eta(E)-x,0\}\Bigr]^{1/\lambda}\right)^{-1}$ \\ & \\

\hline & \\ \textit{Frank t-norms} & $\gamma_E^{F,1}(x) =
\min\Bigl\{\exp(x-\eta(E)),1\Bigr\}$ \\
$T_{\lambda}^{F}$, $\lambda\in ]0,+\infty]$ &
$\gamma_E^{F,+\infty}(x) = \max\Bigl\{\min\{1+x-\eta(E),\, 1\},0\Bigr\}$ \\
& $\gamma_E^{F,\lambda}(x) =
\min\left\{\log_{\lambda}\Bigl(1+(\lambda-1)\exp(x-\eta(E))\Bigr), 1\right\}$ \\ & \\

\hline & \\ \textit{Hamacher t-norms} & $\gamma_E^{H,+\infty}(x) =
\gamma_E^{AA,0}(x)$ \\
$T_{\lambda}^{H}$, $\lambda \in [0,+\infty]$ & $\gamma_E^{H,0}(x)
= \min\left\{\bigl(1 + \eta(E)-x\bigr)^{-1}, 1\right\}$ \\
& $\gamma_E^{H,\lambda}(x) =
\min\left\{\lambda\Bigl(\exp(\eta(E)-x)+\lambda-1\Bigr)^{-1},
1\right\}$ \\
& \\

\hline & \\ \textit{Yager t-norms} & $\gamma_E^{Y,0}(x) =
\gamma_E^{AA,0}(x)$ \\ $T_{\lambda}^{Y}$, $\lambda\in[0,+\infty[$
& $\gamma_{E}^{Y,\lambda}(x) =
\max\left\{\min\left\{1-\Bigl[\max\{\eta(E)-x,0\}\Bigr]^{1/\lambda},\,
1\right\},0\right\}$
\\ & \\

\hline & \\ \textit{Sugeno-Weber t-norms} & $\gamma_E^{SW,-1}(x) =
\gamma_E^{AA,0}(x)$
\\ $T_{\lambda}^{SW}$,
$\lambda\in [-1,+\infty]$ &
$\gamma_E^{SW,0}(x)=\gamma_{E}^{F,+\infty}(x)$ \\ &
$\gamma_E^{SW,+\infty}(x)=\gamma_{E}^{F,1}(x)$ \\ &
$\gamma_E^{SW,\lambda}(x)
=\max\left\{\min\left\{\lambda^{-1}\left((1+\lambda)^{1+x-\eta(E)}-1\right),1\right\},
0\right\}$  \\ & \\
\hline

\end{tabular}\caption{Some well known families of t-norms and their
corresponding parameterized families of
$\tau_{T}$-submeasures}\label{tab_submiery}}
\end{center}
\end{table}

\begin{remark}\rm
Observe that the partial order $\ll$ is a coarsening of the
standard inclusion ordering, i.e., $$\Theta_{L_1,T_1} \ll
\Theta_{L_2,T_2}\,\,\, \Longrightarrow \,\,\, \Theta_{L_1,T_1}
\subset \Theta_{L_2,T_2}.$$ On the other hand, consider for
example $L = K_\infty$. Then $\Theta_{K_\infty, T}$ does not
depend on $T$ (in fact, it consists of probabilistic submeasures
$\gamma$ satisfying $\gamma_E = \gamma_\Omega$ for any non-empty
$E \subset \Omega)$, although there are incomparable t-norms $T_1$
and $T_2$, i.e., $\Theta_{K_\infty,T_1}$ and
$\Theta_{K_\infty,T_2}$ are $\ll$-incomparable.
\end{remark}

\begin{example}\rm
Let $\eta: \Sigma\to\overline{\mathbb{R}}_+$ be a numerical
submeasure on a ring $\Sigma$ of a non-empty set $\Omega$ and
$E\in\Sigma$. Then

(i) $\gamma\in\Theta_{L,M}$, where $L\in\mathcal{L}$, $L\geq K_1$
and $$\gamma_E(x)=\begin{cases}0 & \textrm{for}\,\,x\leq 0, \\
1/2 & \textrm{for}\,\,x\in]0,\eta(E)]; \\ 1 &
\textrm{for}\,\,x>\eta(E),\end{cases}$$

(ii) $\gamma\in\Theta_{L,D}$, where $L\in\mathcal{L}$ is arbitrary
and $$\gamma_E(x)=\frac{x}{x+\eta(E)};$$

(iii) $\gamma\in\Theta_{L,W}$, where $L\in\mathcal{L}$, $L\geq
K_1$ and $$\gamma_E(x)=\max\Bigl\{\min\{1+x-\eta(E),\,
1\},0\Bigr\};$$

(iv) for $$\gamma_E(x) =
\min\left\{\frac{1+x}{1+\eta(E)},1\right\}, \quad x>0,$$ we have
that $\gamma \in \Theta_D$, however it is not an element of
$\Theta_W$, neither $\Theta_\Pi$ nor $\Theta_M$.

More examples of $\tau_T$-submeasures related to some well known
parameterized families of t-norms $T$, see~\cite{KMP}, are
summarized in Table~\ref{tab_submiery}. Note that in all cases we
omit the minimum t-norm $M = T^{AA}_{+\infty} = T^D_{+\infty} =
T^F_0 = T^{SS}_{-\infty} = T^Y_{+\infty}$ (one example of such a
universal submeasure is given in~(i)).
\end{example}

In the context of t-norms (but not limited to this case, as we
will use later) it is very natural to consider the following
simple transformations which often manifest in different applied
fields. Consider the \textit{group $\mathcal{H}$ of automorphisms}
(strictly increasing bijections) of the unit interval $[0,1]$
acting on the class $\mathcal{B}$ of all functions from $[0,1]^2$
to $[0,1]$ as follows
\begin{displaymath}(\Psi_h B)(x,y) =
h^{-1}(B(h(x),h(y))), \quad h\in\mathcal{H},
\end{displaymath} for all $x,y\in[0,1]$. We shall denote
by $\Psi_{\mathcal{H}}$ this class of transformations (an element
of $\Psi_{\mathcal{H}}$ is determined by a function $h \in
\mathcal{H}$). Clearly, $\Psi_\mathcal{H}$ is a group under the
composition with the inverse $\Psi_h^{-1}=\Psi_{h^{-1}}$ and the
identity $\Psi_{\textrm{id}_{[0,1]}}$. The mapping $\Psi:
\mathcal{B}\times\mathcal{H}\to\mathcal{B}$ is the action of the
group $\mathcal{H}$ on $\mathcal{B}$. Since $\Psi_h M=M$ for each
$h\in\mathcal{H}$, then each $\Psi_h$-transform of a universal
submeasure is a universal submeasure as well. Moreover,
$\Theta_{L,M}=\Theta_{L,\Psi_h M}$ for each
$(L,h)\in\mathcal{L}\times\mathcal{H}$. Also it is known,
see~\cite[Proposition 2.6]{DS2}, that the class $\mathcal{T}$ of
t-norms is closed under $\Psi$\footnote{a class $\mathcal{B}$ is
\textit{closed under} $\Psi$, if
$\Psi_h(\mathcal{B})\subset\mathcal{B}$ for each $h\in
\mathcal{H}$}.

\begin{proposition}
Let $h\in\mathcal{H}$. Then
\begin{itemize}
\item[(i)] if $h$ is supermultiplicative, then for each $L_1,
L_2\in\mathcal{L}$ such that $L_1\leq L_2$ holds $\Theta_{L_1,\Pi}
\ll \Theta_{L_2, \Psi_{h}\Pi}$; \item[(ii)] if the function
$1-h(1-x)$ is subadditive, then for each $L_1, L_2\in\mathcal{L}$
such that $L_1\leq L_2$ holds $\Theta_{L_1,W} \ll \Theta_{L_2,
\Psi_{h}W}$; \item[(iii)] if $(L_1,T_1),
(L_2,T_2)\in\mathcal{L}\times \mathcal{T}$ such that
$\Theta_{L_1,T_1} \ll \Theta_{L_2,T_2}$, then
$\Theta_{L_1,\Psi_{h}T_1} \ll \Theta_{L_2,\Psi_{h}T_2}$;
\item[(iv)] for each
$(L,T,h)\in\mathcal{L}\times\mathcal{T}\times\mathcal{H}$ holds
$\Theta_{L,M} \ll \Theta_{L,\Psi_{h}T} \ll \Theta_{L,\Psi_{h}D}$;
\item[(v)] for each $(L,T)\in\mathcal{L}\times\mathcal{T}$ and
each involution $h$ on $[0,1]$ holds: $\gamma\in\Theta_{L,T}$ if
and only if $h\circ\gamma\in\Theta_{L,\Psi_h T}$.
\end{itemize}
\end{proposition}

\begin{example}\rm
Let $\eta$ be a numerical submeasure on $\Sigma$. If
$h(x)=\tan\frac{\pi}{4}x$ for $x\in[0,1]$, then for $L\geq K_1$ we
get $\gamma\in\Theta_{L,\Psi_h W}$, where
$$\gamma_{E}(x)=\max\left\{\min\left\{\frac{4}{\pi}\arctan(1-\eta(E)+x),1\right\},
0\right\}, \quad E\in\Sigma.$$
\end{example}

It is easy to verify that the convex combination of numerical
submeasures is again a numerical submeasure. If we consider the
pseudo-convex combination in the spirit of weighted
quasi-arithmetic mean, the result for probabilistic submeasures
will be the same, i.e., for an arbitrary $L\in\mathcal{L}$ the
weighted quasi-arithmetic mean
$$\mathbf{A}^w_t(x_1,\dots,x_n) = t^{(-1)}\left(\sum_{i=1}^n w_i\,
t(x_i)\right)$$ generated by an additive generator $t$ of a
continuous Archimedean t-norm $T\in\mathcal{T}$ preserves the
class $\Theta_{L,T}$ of probabilistic $\tau_{L,T}$-submeasures.
Here for $i=1,\dots,n$ we consider $x_i\in[0,1]$, $w_i$ are
non-negative weights with $\sum_{i=1}^n w_i=1$ and $t^{(-1)}$ is
the pseudo-inverse function to $t$, see~\cite{KMP} for more
details. Recall that $t:[0,1]\to \overline{\mathbb{R}}_+$ is an
\textit{additive generator} of a continuous Archimedean t-norm $T$
if and only if it is continuous, strictly decreasing and
satisfying $t(1) = 0$. Moreover, its pseudo-inverse $t^{(-1)}:
\overline{\mathbb{R}}_+\to [0,1]$ is given by $$t^{(-1)}(x) =
t^{-1}(\min\{t(0),x\}).$$

\begin{proposition}
Let $L\in\mathcal{L}$ and $t$ be an additive generator of a
continuous Archime\-dean t-norm $T\in\mathcal{T}$. If
$\gamma^{(i)}\in\Theta_{L,T}$ for $i=1,2,\dots,n$, then
$$\gamma = \mathbf{A}^w_t\left(\gamma^{(1)},\dots,\gamma^{(n)}\right)\in\Theta_{L,T}.$$
\end{proposition}

\proof The first two properties~(a') and~(b') of
Definition~\ref{defLT-sub} are easy to verify, therefore we show
only the triangle inequality~(c').

Let $E,F\in\Sigma$. Since $\gamma^{(i)}\in\Theta_{L,T}$, then
$$\gamma_{E\cup F}^{(i)}(L(x,y))\geq
t^{-1}\left(\min\left\{t(0),
t\left(\gamma_E^{(i)}(x)\right)+t\left(\gamma_F^{(i)}(y)\right)\right\}\right),
\quad x,y>0,$$ and we have {\setlength\arraycolsep{2pt}
\begin{eqnarray*}
\gamma_{E\cup F}(L(x,y)) & = & t^{(-1)}\left(\sum_{i=1}^n w_i\,
t\left(\gamma_{E\cup F}^{(i)}(L(x,y))\right)\right) \\ & \geq &
t^{(-1)}\left(\sum_{i=1}^n w_i\, \min\left\{t(0),
t\left(\gamma_E^{(i)}(x)\right)+t\left(\gamma_F^{(i)}(y)\right)\right\}\right)
\\ & \geq & t^{(-1)}\Bigl(t(\gamma_E(x))+t(\gamma_F(y))\Bigr)
\\ & = & T(\gamma_E(x), \gamma_F(y)),
\end{eqnarray*}}thus $\gamma$ is a $\tau_{L,T}$-submeasure on $\Sigma$. \qed

\begin{corollary}
Let $(L,h)\in\mathcal{L}\times\mathcal{H}$ and $t$ be an additive
generator of a continuous Archimedean t-norm $T\in\mathcal{T}$. If
$\gamma^{(i)}\in\Theta_{L,\Psi_h T}$ for $i=1,2,\dots,n$, then
$\gamma=\mathbf{A}^w_t\left(\gamma^{(1)},\dots,\gamma^{(n)}\right)\in\Theta_{L,\Psi_h
T}$.
\end{corollary}

From these observations we state the following open problem:

\begin{problem}
Characterize the class of mappings (aggregation operators)
$\mathcal{M}$ which preserve the class $\Theta_{L,A}$ of
probabilistic submeasures for a fixed
$(L,A)\in\mathcal{L}\times\mathcal{A}$, i.e.,
$\mathcal{M}(\Theta_{L,A})\subseteq \Theta_{L,A}$.
\end{problem}

As it is already known, see~\cite[Theorem 1]{HutMes}, to each
numerical submeasure $\eta$ on $\Sigma$ corresponds
$\gamma\in\Theta_{L,M}$, $L\geq K_1$, in the form
$$\gamma_E(x)=\varepsilon_0(x-\eta(E)), \quad x>0, E\in\Sigma,$$
where the number $\gamma_E(x)$ may be interpreted as the
probability that the value of submeasure $\eta$ of a set
$E\in\Sigma$ is less than $x$. To underline the interesting
relationship between the probabilistic $\tau_{L,T}$-submeasure
$\gamma$ and the numerical submeasure $\eta$ on $\Sigma$ we give
the following result which improves and generalizes~\cite[Theorem
4]{HutMes}. For the sake of completeness we give its short direct
proof here.

\begin{theorem}\label{thmnumerical}
Let $L\leq K_1$ and $\gamma\in \Theta_{L,T_1}$. If $t$ is an
additive generator of a continuous Archimedean t-norm $T$ such
that $T\leq T_1$, then a mapping $\eta_{\gamma, t}: \Sigma \to
\mathbb{R}_+$ given by
$$\eta_{\gamma, t}(E) = \sup\{z \in \mathbb{R}_+; \, t(\gamma_{E}(z)) \geq
z\}$$ is a numerical submeasure.
\end{theorem}

\proof The equality $\eta_{\gamma, t}(\emptyset) = 0$ and the
monotonicity of $\eta_{\gamma, t}$ are obvious. Moreover, it is
evident that $\gamma$ is an element of $\Theta_{K_1,T}$, and hence
for $E,F\in\Sigma$ we have {\setlength\arraycolsep{2pt}
\begin{eqnarray*}
& \phantom{=} & \eta_{\gamma, t}(E\cup F) = \sup\{z\in
\mathbb{R}_+;\, t(\gamma_{E\cup F}(z)) \geq z\} \\ & \leq &
\sup\Bigl\{z\in \mathbb{R}_+;
\,t(T(\gamma_E(x),\gamma_F(z-x)))\geq x+z-x \,\, \textrm{for
some}\,\, x \in [0,z]\Bigr\} \\ & = & \sup\Bigl\{z\in
\mathbb{R}_+;\, \min\{t(0), t(\gamma_E(x)) + t(\gamma_F(z-x))\}
\geq z \,\,\textrm{for some}\,\, x \in [0,z]\Bigr\} \\ & \leq &
\eta_{\gamma, t}(E) + \eta_{\gamma, t}(F),
\end{eqnarray*}}which proves that $\eta_{\gamma,t}$ is a numerical submeasure on $\Sigma$. \qed

Now we will consider the Fr\'echet-Nikodym topology
$\Gamma(\gamma)$ generated by probabilistic submeasure $\gamma$ on
$\Sigma$. This notion was introduced and studied by Drewnowski
in~\cite{Drewnowski} for numerical submeasures on a ring of sets.
Recall that a topology $\sigma$ on a ring $\Sigma$ is said to be a
\textit{ring topology} if the mappings $(E,F)\to E\triangle F$ and
$(E,F)\to E\cap F$ of $\Sigma\times\Sigma\to\Sigma$ are continuous
(with respect to the product topology on $\Sigma\times\Sigma$). A
ring topology $\sigma$ is said to be a \textit{Fr\'echet-Nikodym
topology} on $\Sigma$ if for each $\sigma$-neighborhood $U$ of
$\emptyset$ in $\Sigma$ there is a $\sigma$-neighborhood $V$ of
$\emptyset$ in $\Sigma$ such that $F\subset U$ for all $F\subseteq
E\in V$, $F\in\Sigma$. In particular, a family $\{\eta_i; i\in
I\}$ of numerical submeasures on $\Sigma$ defines a
Fr\'echet-Nikodym topology $\Gamma(\eta_i; i\in I)$ and
conversely, for each Fr\'echet-Nikodym topology $\Gamma$ on
$\Sigma$ there is a family $\{\zeta_j; j\in J\}$ of numerical
submeasures on $\Sigma$ such that $\Gamma=\Gamma(\zeta_j; j\in
J)$.

Define the set function $\rho: \Sigma\times\Sigma \to \Delta^+$ by
$\rho(E,F)=\gamma_{E\triangle F}$ where $\gamma\in\Theta_{L,T}$.
Then {\setlength\arraycolsep{2pt}
\begin{eqnarray*}
\rho_{E,F}(L(x,y)) & = & \gamma_{E\triangle F}(L(x,y)) \geq
\gamma_{(E\triangle G)\cup (G\triangle F)}(L(x,y)) \nonumber\\ &
\geq & T(\gamma_{E\triangle G}(x), \gamma_{G\triangle F}(y)) =
T(\rho_{E,G}(x), \rho_{G,F}(y)),
\end{eqnarray*}}which means, in the other words, that $\rho$ is an
$L$-Menger pseudo-metric on $\Sigma$. Moreover, $\rho$ is
translation invariant, i.e., $$\rho_{E,F} = \gamma_{E\triangle F}
= \gamma_{(E\triangle G)\triangle (G\triangle F)} =
\rho_{E\triangle G, G\triangle F}.$$ Thus, the triple $(\Sigma,
\rho, \tau_{L,T})$ is an $L$-Menger probabilistic pseudo-metric
space, see~\cite[Theorem 3.2]{HalHutMes}. Since for $E_1, E_2,
F_1, F_2\in\Sigma$ holds
$$(E_1\cap F_1)\triangle (E_2\cap F_2) \subset (E_1\triangle E_2)
\cup (F_1\triangle F_2),$$ then we get
{\setlength\arraycolsep{2pt}
\begin{eqnarray}\label{(1)}
\rho_{E_1\cap F_1, E_2\cap F_2}(x) & = & \gamma_{(E_1\cap
F_1)\triangle (E_2\cap F_2)}(x) \geq \gamma_{(E_1\triangle E_2)
\cup (F_1\triangle F_2)}(L(z,z)) \\ & \geq &
T(\gamma_{E_1\triangle E_2}(z), \gamma_{F_1\triangle F_2}(z)) =
T(\rho_{E_1,E_2}(z), \rho_{F_1,F_2}(z)), \nonumber
\end{eqnarray}}where $L(z,z)<x$. If $(E_n,F_n)\to (E,F)$ in topology $\Gamma(\gamma)$,
then $E_n\to E$ and $F_n\to F$. Thus, $\rho_{E_n, E}(z)\to 1$ and
$\rho_{F_n, F}(z)\to 1$. Moreover, if we consider a continuous
t-norm $T$, then from~(\ref{(1)}) we get $\rho_{E_n\cap E, F_n\cap
F}(z)\to 1$ for each $x>0$. In fact, it proves continuity of
$\cap$ in the product topology $\Sigma\times\Sigma$. These
observations lead to the following result.

\begin{proposition}
Let $(L,T)\in\mathcal{L}\times\mathcal{T}$, where $T$ is a
continuous t-norm and $\gamma\in\Theta_{L,T}$. For $\varepsilon>0$
and $\delta>0$ put
$$\mathcal{B}(\varepsilon,\delta)=\{E\in\Sigma;\,\,\,
\gamma_E(\varepsilon)>1-\delta\}.$$ Then

\begin{itemize}
\item[(i)] $\mathfrak{B} = \{\mathcal{B}(\varepsilon, \delta);\,\,
\varepsilon>0, \delta>0\}$ is a normal base of neighborhoods of
$\emptyset$ for the Fr\'echet-Nikodym topology $\Gamma(\gamma)$;
\item[(ii)] $(\Sigma, \triangle, \cap, \Gamma(\gamma))$ is a
topological ring of sets. 
\end{itemize}
\end{proposition}

\section{Semi-copula-based submeasures}\label{sectioncopulas}

In what follows we consider the natural extension/modification of
t-norms: copulas, quasi-copulas and semi-copulas, see~\cite{DS}.
Recall that a \textit{semi-copula} is an aggregation function $S:
[0,1]^2 \to [0,1]$ with 1 as its neutral element. Denote by
$\mathcal{S}$ the set of all semi-copulas and $\mathcal{S}_c$ the
set of all continuous semi-copulas. A \textit{quasi-copula} $Q$ is
a 1-Lipschitz semi-copula, i.e., a semi-copula $Q$ satisfying
$$|Q(x,y)-Q(x',y')|\leq |x-x'|+|y-y'|$$ for all $x,x',y,y'\in[0,1]$. The set of all
quasi-copulas will be denoted by $\mathcal{Q}$. A semi-copula $C$
which is 2-increasing, i.e., for each $x,y,x',y'\in ]0,1]$ such
that $x\leq x'$ and $y\leq y'$ holds
$$C(x',y')-C(x,y')-C(x',y)+C(x,y)\geq 0,$$ is called a \textit{copula}. Denote by $\mathcal{C}$
the set of all copulas. Then $\mathcal{C}\subset \mathcal{Q}
\subset \mathcal{S}$. 
All these sets are partially ordered equipped with the usual
point-wise order $\leq$ between real functions. Clearly, for each
$h\in\mathcal{H}$ the mapping $\Psi_h$ is order-preserving on
$\mathcal{S}$ and for a given $h\in\mathcal{H}$ the partially
ordered set
$$\mathcal{K}=(\{S\in\mathcal{S}; \Psi_h S = S\}, \leq)$$ is a
complete lattice (by Knaster-Tarski theorem). As we have already
mentioned, $M\in\mathcal{K}$, but also $D\in\mathcal{K}$.



For the class of Archimedean copulas, i.e., copulas of the form
$$C(x,y)=\varphi^{[-1]}(\varphi(x)+\varphi(y))$$ for all
$x,y\in[0,1]$ where $\varphi:[0,1]\to[-\infty,+\infty]$ is a
continuous, strictly decreasing convex function with
$\varphi(1)=0$ and the pseudo-inverse $\varphi^{[-1]}$ (such a
function is called an \textit{additive generator} of $C$,
cf.~\cite{Nelsen}) we immediately have the following
characterization.

\begin{proposition}
Let $\eta$ be a numerical submeasure on $\Sigma$. If $\varphi$ is
an additive generator of $C\in\mathcal{C}$, then
$\gamma\in\Theta_{C}$, where
$$\gamma_E(x)=\varphi^{[-1]}(\eta(E)-x), \quad E\in\Sigma.$$ Moreover, for each $h\in\mathcal{H}$ holds
$\gamma\in\Theta_{\Psi_h C}$, where $$\gamma_E(x)=(\varphi\circ
h)^{[-1]}(\eta(E)-x), \quad E\in\Sigma.$$
\end{proposition}

Easily it is possible to state the analogous result for the
multiplicative generator of $C\in\mathcal{C}$.

\begin{example}\rm
Let $\eta$ be a numerical submeasure on $\Sigma$ and $E\in\Sigma$.
Then

(i) $\gamma\in\Theta_{C^{GH}_\lambda}$, where
$$\gamma_{E}(x)=\exp\left(-\Bigl[\max\{\eta(E)-x,0\}\Bigr]^{1/\lambda}\right)$$
corresponds to the \textit{Gumbel-Hougaard family} of (strict)
copulas $C^{GH}_\lambda$ given by
$$C^{GH}_\lambda(u,v)=\exp\left(-\left[(-\ln u)^\lambda + (-\ln v)^\lambda\right]^{1/\lambda}\right),$$
with $\lambda\in [1,+\infty[$, see~\cite{Nelsen}; for $\lambda=1$
we have the independence copula $\Pi$ (and the corresponding
$\gamma\in\Theta_\Pi$) -- clearly,
$\Theta_\Pi\in\Theta_\mathcal{C}$ and for each $h\in\mathcal{H}$
we have $\Theta_{\Psi_h \Pi}\in\Theta_\mathcal{T}$;

(ii) $\gamma\in\Theta_{C_\lambda}$, where
$$\gamma_{E}(x)=\max\left\{\min\left\{\frac{1-\eta(E)+x}{1+(\lambda-1)(\eta(E)-x)},1\right\},0\right\}$$
corresponds to the family of (non-strict) copulas
\begin{displaymath}
C_\lambda(u,v) = \max\left\{\frac{\lambda^2
uv-(1-u)(1-v)}{\lambda^2-(\lambda-1)^2(1-u)(1-v)}, 0\right\},
\quad \lambda\in [1,+\infty[.
\end{displaymath}
\end{example}

Observe that similarly as in the case of t-norms $\mathcal{T}$,
the weighted quasi-arithmetic mean
$$\mathbf{A}^w_\varphi(x_1,\dots,x_n) = \varphi^{[-1]}\left(\sum_{i=1}^n w_i\,
\varphi(x_i)\right)$$ generated by an additive generator $\varphi$
of an Archimedean copula $C\in\mathcal{C}$ preserves the class
$\Theta_{C}$ of probabilistic copula-based submeasures (even their
generalization involving an arbitrary $L\in\mathcal{L}$).

Given a class $\mathcal{B}$ of functions from $[0,1]^2$ to
$[0,1]$, we shall denote by $\Psi_{\mathcal{H}}(\mathcal{B})$ the
class of operators obtained by transforming all elements of
$\mathcal{B}$ by all elements of $\Psi_{\mathcal{H}}$. Since the
classes $\mathcal{S}$ and $\mathcal{S}_c$ are closed under $\Psi$,
cf.~\cite{DS2}, we have the following relations
$$\Psi_{\mathcal{H}}(\mathcal{C}) \subset \Psi_{\mathcal{H}}(\mathcal{Q}) \subset
\Psi_{\mathcal{H}}(\mathcal{S}_c)=\mathcal{S}_c \subset
\mathcal{S} = \Psi_{\mathcal{H}}(\mathcal{S}).$$ Recall that the
identity $I\in\Psi_\mathcal{H}$. Moreover, $\mathcal{C}\subset
\Psi_{\mathcal{H}}(\mathcal{C})$ and $\mathcal{Q}\subset
\Psi_{\mathcal{H}}(\mathcal{Q})$, see~\cite{APS}.

\begin{proposition}
If $S_1, S_2\in\mathcal{S}$ such that
$\Theta_{S_1}\ll\Theta_{S_2}$, then for each $h\in\mathcal{H}$ it
holds $\Theta_{\Psi_{h}S_1} \ll \Theta_{\Psi_{h}S_2}$. Moreover,
$\Theta_{\mathcal{S}}=\Theta_{\Psi_{\mathcal{H}}(\mathcal{S})}$
and
$\Theta_{\mathcal{S}_c}=\Theta_{\Psi_{\mathcal{H}}(\mathcal{S}_c)}$.
\end{proposition}

In what follows we are interested in lattice structure of
submeasure spaces in $\Theta_\mathcal{S}$. As shown in~\cite{DMP},
the class $\mathcal{S}$ of semi-copulas constitutes the lattice
completion of the class $\mathcal{T}$ of t-norms, in the sense
that every semi-copula may be represented as the point-wise
supremum and infimum of a suitable subset of t-norms. Let $\vee$
and $\wedge$ denote the point-wise supremum and infimum,
respectively. Observe that if $\gamma$ is a $\tau_{S_1}$- and
$\tau_{S_2}$-submeasure for some $S_1, S_2\in \mathcal{S}$, then
$\gamma$ is a $\tau_{S_1\vee S_2}$- as well as $\tau_{S_1 \wedge
S_2}$-submeasure. Thus, for $S_1, S_2\in\mathcal{S}$ put
$$\Theta_{S_1} \sqcup \Theta_{S_2} = \Theta_{S_1 \wedge S_2}\quad
\textrm{and} \quad \Theta_{S_1} \sqcap \Theta_{S_2} = \Theta_{S_1
\vee S_2}.$$ It is easy to see that $\sqcup$ and $\sqcap$ are
\textit{lattice operations}. Since $(\mathcal{S}, \leq, \vee,
\wedge)$ is a complete lattice, see~\cite{DMP}, then we have the
following observation.

\begin{proposition}
The family $\Theta_{\mathcal{S}}$ of all probabilistic submeasure
spaces is a distributive lattice.
\end{proposition}

\proof Indeed, for $S_1, S_2, S_3 \in \mathcal{S}$, we have:
{\setlength\arraycolsep{2pt}
\begin{eqnarray*}
\Theta_{S_1} \sqcup (\Theta_{S_2} \sqcap \Theta_{S_3}) & = &
\Theta_{S_1} \sqcup \Theta_{S_2 \vee S_3} = \Theta_{S_1 \wedge
(S_2 \vee S_3)} = \Theta_{(S_1 \wedge S_2) \vee (S_1 \wedge S_3)}
\\ & = &  \Theta_{S_1 \wedge S_2} \sqcap \Theta_{S_1 \wedge S_3} =
(\Theta_{S_1} \sqcup \Theta_{S_2}) \sqcap (\Theta_{S_1} \sqcup
\Theta_{S_3}).
\end{eqnarray*}}Analogously for $\Theta_{S_1} \sqcap (\Theta_{S_2} \sqcup \Theta_{S_3})$.
By~\cite[Theorem~2.2]{Birkhoff} $\Theta_{\mathcal{S}}$ is a
distributive lattice. \qed

Since for each $S\in\mathcal{S}$ holds $\Theta_{M}\ll \Theta_{S}
\ll \Theta_{D}$, then $\Theta_{M}$ is bottom and $\Theta_{D}$ is
top in the lattice $\Theta_{\mathcal{S}}$, thus
$\Theta_{\mathcal{S}}$ is a \textit{bounded distributive lattice}.

\begin{proposition}\label{Theorem 2.4}
For every $S_1 \in \mathcal{S}$, the set $$\mathfrak{I}_{S_1} =
\{\Theta_{S} \in \Theta_{\mathcal{S}};\ \Theta_{S} \ll
\Theta_{S_1},\, S \in \mathcal{S}\}$$ is an ideal in
$\Theta_{\mathcal{S}}$.
\end{proposition}


\proof First, observe that for $S_1\in\mathcal{S}$ the set
$\mathfrak{I}_{S_1}$ is the set of all $\tau_S$-submeasures
related to a semi-copula $S$ which are also
$\tau_{S_1}$-submeasures.

Let $\Theta_{S_2} \in \mathfrak{I}_{S_1}$, i.e., $\Theta_{S_2} \ll
\Theta_{S_1}$ and let $\Theta_{S_3} \ll \Theta_{S_2}$. From it
follows that $S_1\leq S_2$ and $S_2\leq S_3$. Thus, $S_1\leq S_3$
which shows that $\Theta_{S_3} \ll \Theta_{S_1}$, i.e.,
$\Theta_{S_3} \in \mathfrak{I}_{S_1}$.

Let $\Theta_{S_2}, \Theta_{S_3} \in \mathfrak{I}_{S_1}$, i.e.,
$S_1\leq S_2$ and $S_1\leq S_3$. Since $S_1\leq S_2 \wedge S_3$,
then $\Theta_{S_2} \sqcup \Theta_{S_3} = \Theta_{S_2 \wedge S_3}
\ll \Theta_{S_1}$, i.e., $\Theta_{S_2} \sqcup
\Theta_{S_3}\in\mathfrak{I}_{S_1}$. Therefore,
$\mathfrak{I}_{S_1}$ is an ideal in $\Theta_{\mathcal{S}}$. \qed

Dually to Theorem~\ref{Theorem 2.4}, we obtain the following
corollary.


\begin{corollary}\label{Corollary 2.5}\rm
For every $S_2 \in \mathcal{S}$, the set $$\mathfrak{F}_{S_2} =
\{\Theta_{S} \in \Theta_{\mathcal{S}};\ \Theta_{S_2} \ll
\Theta_{S}, \, S \in \mathcal{S}\}$$ is a filter in
$\Theta_{\mathcal{S}}$.
\end{corollary}

From it follows that for $S_1,S_2\in\mathcal{S}$ such that
$S_1\leq S_2$ the set $$[\Theta_{S_2}, \Theta_{S_1}] =
\mathfrak{I}_{S_1}\cap \mathfrak{F}_{S_2}$$ is a order interval in
$\Theta_\mathcal{S}$. 

\begin{theorem}\label{proplatticeThetaSQ}
$(\Theta_{\mathcal{S}}, \ll, \sqcap, \sqcup, \Theta_D, \Theta_M)$
is a complete lattice.
\end{theorem}

\proof Let $\Theta_{\mathcal{P}}$ be any subset of
$\Theta_{\mathcal{S}}$ and put $\sqcup \Theta_{\mathcal{P}} =
\Theta_{\wedge\mathcal{P}}$, $\sqcap \Theta_{\mathcal{P}} =
\Theta_{\vee\mathcal{P}}$, where $$\vee\mathcal{P}(x,y) =
\sup\{P(x,y); P\in\mathcal{P}\} \quad \wedge\mathcal{P}(x,y) =
\inf\{P(x,y); P\in\mathcal{P}\}$$ for each $(x,y)\in[0,1]^2$.
Since $(\mathcal{S}, \leq, \vee, \wedge)$ is complete, then for
each $\mathcal{P}\subseteq\mathcal{S}$ holds $\vee
\mathcal{P}\in\mathcal{S}$ and $\wedge\mathcal{P}\in\mathcal{S}$.
Thus $\Theta_{\vee \mathcal{P}}\in\Theta_{\mathcal{S}}$ and
$\Theta_{\wedge \mathcal{P}}\in\Theta_{\mathcal{S}}$. \qed

\begin{remark}\rm
Since $(\mathcal{Q}, \leq, \vee, \wedge)$ is a complete lattice,
see~\cite{NU-F}, all the above assertions hold also for
$\Theta_{\mathcal{Q}}$, i.e., $(\Theta_\mathcal{Q}, \ll, \sqcap,
\sqcup, \Theta_W, \Theta_M)$ is a complete sublattice of
$(\Theta_{\mathcal{S}}, \ll, \sqcap, \sqcup, \Theta_D, \Theta_M)$,
where for each $Q_1\in\mathcal{Q}$ the set $$\mathfrak{I}_{Q_1} =
\{\Theta_{Q} \in \Theta_{\mathcal{Q}};\ \Theta_{Q} \ll
\Theta_{Q_1},\, Q \in \mathcal{Q}\}$$ is an ideal in
$\Theta_{\mathcal{Q}}$ and for each $Q_2 \in \mathcal{Q}$, the set
$$\mathfrak{F}_{Q_2} = \{\Theta_{Q} \in \Theta_{\mathcal{Q}};\
\Theta_{Q_2} \ll \Theta_{Q}, \, Q \in \mathcal{Q}\}$$ is a filter
in $\Theta_{\mathcal{Q}}$.
\end{remark}

\section*{Acknowledgement}
This work was partially supported by research grants APVV-0073-10
and P402/11/0378. The second author has been on a postdoctoral
stay at the Departamento de Matem\'aticas, CINVESTAV del IPN
(M\'exico) when investigating some parts of the work presented
therein and finishing the paper. He therefore gratefully
acknowledges the hospitality and support of the mathematics
department of CINVESTAV on this occasion and partial support of
the Research Project VVGS 45/10-11.

\vspace{5mm}

\noindent \small{\textsc{Lenka Hal\v cinov\'a} \\ Institute of
Mathematics, Faculty of Science, Pavol Jozef \v Saf\'arik
University in Ko\v sice, {\it Current address:} Jesenn\'a 5,
040~01 Ko\v sice, Slovakia \\ {\it E-mail address:}
lenka.halcinova@student.upjs.sk}

\vspace{5mm}

\noindent \small{\textsc{Ondrej Hutn\'ik} \\ Institute of
Mathematics, Faculty of Science, Pavol Jozef \v Saf\'arik
University in Ko\v sice, Jesenn\'a 5, 040~01 Ko\v sice, Slovakia
\\ {\it E-mail address:} ondrej.hutnik@upjs.sk} \\ AND \\ Departamento de
Matem\'aticas, CINVESTAV del IPN, \\ {\it Current address:}
Apartado Postal 14-740, 07000, M\'exico, D.F., M\'exico
\\ {\it E-mail address:} hutnik@math.cinvestav.mx

\vspace{5mm}

\noindent \small{\textsc{Radko Mesiar} \\ Department of
Mathematics and Descriptive Geometry, Faculty of Civil
Engineering, Slovak University of Technology, \\ {\it Current
address:} Radlinsk\'eho 11, 813~68 Bratislava, Slovakia \\ AND
\\ Institute of Information Theory and Automation, Academy of
Sciences of the Czech Republic, P.O.Box 18, Prague 18208, Czech
Republic. \\ {\it E-mail address:} mesiar@math.sk}

\end{document}